\documentclass{amsart}

\usepackage{amsmath,amsthm,amsfonts,amscd,amssymb} 
\usepackage{verbatim}      	% Allows quoting source with commands.
\usepackage{epsfig}         	% Allows inclusion of eps files.
\usepackage{url}		% Allows good typesetting of web URLs.
\usepackage{epic,eepic,latexsym}
\usepackage{graphicx}
\usepackage{tikz-cd}
\usepackage{color}
\usepackage{lcd}
\usepackage{mathrsfs}
\usepackage[centering,text={15.5cm,22cm},
        marginparwidth=20mm,dvips]{geometry}
\usepackage[linktocpage]{hyperref}
\usepackage{lscape}
\usepackage{tabularx}
\usepackage{marginnote}
\usepackage[all]{xy}
\usepackage{enumitem}
\usepackage{stmaryrd}

\DeclareMathOperator{\Spec}{Spec}

\DeclareMathOperator{\prin}{prin}

\DeclareMathOperator{\Spf}{Spf}
\DeclareMathOperator{\Supp}{Supp}

\DeclareMathOperator{\up}{up}

\let\?=\overline
\let\:=\colon

\let\bb=\mathbb

\let\f=\mathfrak
\let\s=\mathcal

\let\wh=\widehat

\newcommand {\kk} {\Bbbk}

\theoremstyle{plain}% default
 \newtheorem{thm}{Theorem}

\theoremstyle{definition}

 \newtheorem{egs}[thm]{Examples}
 
\theoremstyle{remark}

\title{Theta bases are atomic}
\author{Travis Mandel}
\address{University of Utah\\
Department of Mathematics\\
155 S 1400 E RM 233\\
Salt Lake City, UT, 84112-0090}
\email{mandel{\char'100}math.utah.edu}
\thanks{The author was supported by the Center of Excellence Grant ``Centre for Quantum Geometry of Moduli Spaces'' from the Danish National Research Foundation (DNRF95) and later by the National Science Foundation RTG Grant DMS-1246989.}

\begin{document}

\begin{abstract}
Fock and Goncharov conjectured that the indecomposable universally positive (i.e., atomic) elements of a cluster algebra should form a basis for the algebra.  This was shown to be false by Lee-Li-Zelevinsky.  However, we find that the theta bases of Gross-Hacking-Keel-Kontsevich do satisfy this conjecture for a slightly modified definition of universal positivity in which one replaces the positive atlas consisting of the clusters by an enlargement we call the scattering atlas.  In particular, this uniquely characterizes the theta functions.
\end{abstract}

\maketitle

\setcounter{tocdepth}{0}
%\tableofcontents  

\section*{}
A cluster variety $U$, as defined in \cite{FG1}, is a scheme constructed by gluing together a collection of algebraic tori $\Spec \kk[M]$, called {\bf clusters}, via certain birational automorphisms called mutations.  A nonzero $f\in \Gamma(U,\s{O}_U)$ is called {\bf universally positive} if its restriction to each cluster is a Laurent polynomial $\sum_{m\in M} a_m z^m$ with non-negative integer coefficients, and {\bf indecomposable} or {\bf atomic} if it cannot be written as a sum of two other universally positive functions.

Fock and Goncharov predicted \cite[Conjecture 4.1]{FG1} that the atomic functions on $U$ form an additive basis for $\Gamma(U,\s{O}_U)$.  However, \cite{LLZ} showed this to be false by showing that the atomic functions are often linearly dependent, even in rank $2$.   Nevertheless, \cite{GHKK} constructed a canonical topological basis of ``theta functions'' $\{\vartheta_m\}_{m\in M}$ for a topological algebra $A$ that should be viewed as a formal version of (a base extension of) $\Gamma(U,\s{O}_U)$---e.g., $A=\wh{\up(\?{\s{A}}_{\prin}^s)} \otimes_{\kk[N_s^+]} \kk[N]$ as in \cite[\S 6]{GHKK} when $U$ is the $\s{A}$-variety, or $A=\wh{B}$ as in \cite[\S 2.4]{Man3}.  These {\bf theta bases} satisfy many of the properties predicted by \cite{FG1} (e.g., universal positivity, being parametrized by $M$), and in many cases (when ``the full Fock-Goncharov conjecture holds''), they extend to bases of $\Gamma(U,\s{O}_U)$.

The construction of $\{\vartheta_m\}_{m\in M}$ (cf. \cite{GHKK} for the details) involves a ``scattering diagram'' $\f{D}$ in $M_{\bb{R}}:=M\otimes \bb{R}$.  $\f{D}$ consists of a set of codimension $1$ ``walls'' in $M_{\bb{R}}$ (with attached functions), the union of which forms the support of $\f{D}$, denoted $\Supp(\f{D})$.  Each $Q\in M_{\bb{R}}\setminus \Supp(\f{D})$ determines an inclusion $\iota_Q$ of $A$ into a certain Laurent series ring $\wh{R[M]}$---a localization of a completion of a base extension of $\kk[M]$, or more precisely,  $\wh{R[M]}=\wh{\kk[\sigma]}\otimes_{\kk[\sigma]} \kk[M']$, where $M'$ is some lattice containing $M$ (e.g., $M\oplus M^*$ for the cluster algebra with principal coefficents), $\sigma$ is some cone in $M'$, and $\wh{\kk[\sigma]}$ denotes the completion of $\kk[\sigma]$ with respect to its maximal ideal.  We say a nonzero $f\in A$ is {\bf universally positive with respect to the scattering atlas} if for every $Q\in M_{\bb{R}}\setminus \Supp(\f{D})$, $\iota_Q(f)\in \wh{R[M]}$ is a formal Laurent series with non-negative integer coefficients.  Such an $f$ is called {\bf atomic with respect to the scattering atlas} if it cannot be written as a sum of two other elements which are universally positive with respect to the scattering atlas.

\begin{thm}\label{1}
The theta functions are exactly the atomic elements of $A$ with respect to the scattering atlas.
\end{thm}

To justify the terminology, recall that for $Q_1,Q_2\in M_{\bb{R}}\setminus \Supp(\f{D})$, the $\iota_{Q_i}$'s are related by a ``path-ordered product'' (cf. \cite[\S 1.1]{GHKK}).  These generalizations of mutations are automorphisms of $\wh{R[M]}$ that take a monomial $z^m$ to a formal positive rational function, by which we mean a quotient $f_1/f_2$ for two nonzero Laurent series $f_1,f_2 \in \wh{R[M]}$ with non-negative integer coefficients.  Hence, the charts $\Spf \wh{R[M]} \hookrightarrow \Spf A$ induced by the $\iota_Q$'s form a {\bf formal positive atlas} on $\Spf A$---i.e., a positive atlas in the sense of \cite[Def. 1.1]{FG1}, except that the split algebraic tori there are replaced with the copies of the formal algebraic torus $\Spf \wh{R[M]}$, and the transition maps now preserve {\it formal} positive rational functions.  This formal positive atlas is what we call the {\bf scattering atlas}.

We note that when $U$ is a cluster $\s{A}$-variety, there is a subset of the chambers of $M_{\bb{R}}\setminus \Supp(\f{D})$, called the cluster complex, such that the corresponding charts of the scattering atlas are just the clusters (restricted to $\Spf A$).  The scattering atlas may therefore be viewed as an enlargement of the positive atlas consisting of the clusters, which we call the {\bf cluster atlas}.  Note that the cluster atlas is the one considered by \cite{FG1}.  It is not clear in general whether the theta functions are indecomposable with respect to the cluster atlas, but \cite{CGMMRSW} has shown this for rank $2$.

\begin{egs}
For cluster algebras of finite or affine type, the cluster complex is dense in $M_{\bb{R}}$, so the scattering atlas and cluster atlas agree in these cases.  Similarly, it is known that the cluster $\s{A}$-variety associated to the Markov quiver $
\begin{tikzcd}[arrow style=tikz,>=stealth,column sep=.04em,row sep=.70 em,nodes={scale=0.5}]
& \bullet \arrow[dl,shift left=.3ex]
  \arrow[dl,shift right=.3ex] \\
\bullet \arrow[rr,shift left=.3ex]\arrow[rr,shift right=.3ex]
&& \bullet   \arrow[ul,shift left=.3ex]
  \arrow[ul,shift right=.3ex]
\end{tikzcd}
$ admits two cluster structures (to be described in detail in \cite{Zhou}), and the union of the corresponding cluster complexes is dense in $M_{\bb{R}}$ (these are the cones over the hemispheres $S^+$ and $S^-$ described in \cite[\S 2.2]{FG2}).  Hence, being universally positive with respect to the scattering atlas here is equivalent to being universally positive with respect to both cluster atlases.
\end{egs}

\begin{proof}[Proof of Theorem \ref{1}]
The fact that the theta functions are universally positive with respect to the scattering atlas was already observed in \cite{GHKK}.  It is an easy consequence of the positivity of the scattering diagram in their Theorem 1.28.  To show indecomposability, it thus suffices to show that for any $f\in A$ universally positive with respect to the scattering atlas, the expansion $f=\sum_{m\in M} a_m \vartheta_m$ has non-negative integer coefficients. For $Q$ sufficiently close to $p\in M$ and $\iota_Q(f)=\sum_{m\in M} c_m z^m\in \wh{R[M]}$, the proof of \cite[Prop. 6.4]{GHKK} shows that $a_p=c_p$.  This is indeed a non-negative integer by the positivity assumption on $f$.
\end{proof}

The same argument proves that the {\it quantum} theta bases (as constructed in \cite{Man3}) are exactly the atomic elements of the corresponding quantum algebras with respect to the scattering atlas, assuming universal positivity of these bases.  This positivity fails in general but is expected for skew-symmetric cases, where positivity of the cluster variables in all clusters was recently proved in \cite{Dav}. 

\subsection*{Acknowledgements}
I would like to thank Li Li for encouraging me to investigate the indecomposability of the theta functions, as well as Sean Keel for his advice on writing this paper.

\bibliographystyle{amsalpha}  % Here the bibliography 		     %
\bibliography{mandel}        % is inserted.			     %
\index{Bibliography@\emph{Bibliography}}%

\end{document}